\documentclass[10pt]{article}
\newtheorem{theorem}{THEOREM}
\newtheorem{proposition}{PROPOSITION}
\newtheorem{corollary}{COROLLARY}
\begin{document}

\renewcommand{\theequation}{\thesection.\arabic{equation}}

\title{Controllability properties for the one-dimensional Heat equation
under multiplicative or nonnegative additive controls \\ with local
mobile support}
\author{Luis A. Fern\'andez$^1$\thanks{\ The work of the first author was partially
supported by the Spanish Ministry of Science and Innovation under grant $MTM2008-04206$.}
\ and Alexander Y. Khapalov$^2$ \\
$^1$Dep. de Matem\'aticas, Estad\'{\i}stica y
Computaci\'on, \\ Avda de los Castros, s/n,\\
 Universidad de Cantabria,  39005-Santander, SPAIN\\ $^2$Department of Mathematics, \\
Washington State University,  Pullman, WA 99164-3113, USA\\ E-mail:
lafernandez@unican.es, khapala@math.wsu.edu}
\date{}
\maketitle

\begin{abstract}
We discuss several new results on nonnegative approximate
controllability for the one-dimensional Heat equation governed by
either multiplicative or nonnegative additive control, acting within
a proper subset of the space domain at every moment of
time. Our methods allow us to link these two types of controls to some extend. The main results include approximate
controllability properties both for the static and mobile control
supports.
\end{abstract}


{\bf AMS(MOS) subject classifications.} 35K05, 35K20, 93B05.
\bigskip

\maketitle
\section{Problem setting}
\setcounter{equation}{0} In this paper we study the approximate
controllability properties of the one dimensional Heat equation
governed either by an additive nonnegative locally distributed
control or by a multiplicative control of local support. For
simplicity, we will assume that the space domain is the interval
$(0,1)$. More precisely, we consider the following two Dirichlet
problems:

$$
\left\{
\begin{array}{lcll}
y_{t}(x,t) &  =   & y_{xx}(x,t) + u(x,t)\chi_{\omega(t)}(x), & {\rm
in} \; Q_T = (0,1)  \times (0,T),
\\y(0,t) & = & y(1,t) =  0,   &  {\rm in } \;\; (0,T),
\\
y(x,0) &  = & y_0(x),  & {\rm in } \;\; (0,1),
\end{array}
\right. \eqno(S1)$$ and
$$ \left\{
\begin{array}{lcll}
y_{t}(x,t) &  =   & y_{xx}(x,t) + v(x,t) y(x,t) \chi_{\omega(t)}(x),
& {\rm in} \; Q_T,
\\y(0,t) & = & y(1,t) =  0,   &  {\rm in } \;\; (0,T),
\\
y(x,0) &  = & y_0 (x),  & {\rm in } \;\; (0,1),
\end{array}
\right. \eqno(S2)$$ where $y_0 \in L_+^2 (0,1) = \{\phi \in L^2
(0,1) : \phi(x) \geq 0 \;\; {\rm a.e. \;\; in } \; (0,1)\}$, with $y_0
\not\equiv 0$ for $(S2)$. In the above, $u$ is  the additive locally
distributed control, while $v$ is the multiplicative (or {\em
bilinear}) control (in their respective systems), each one supported
at each time $ \; t \; $ in a proper subinterval $\omega (t)$ of $
[0,1]$. We will assume that
\begin{equation} \omega(t) = (r(t), r(t) + l), \label{1.1}\end{equation} where the length $ l \in (0,
1) $ is fixed and $r: [0,T] \longrightarrow [0, 1-l]$ is a function
that defines the position of the support at time $t$. We call the
support $ \omega(t) $ {\em static} if it does not change in time
and {\em mobile} if it does. The symbol $ \chi_{\omega(t)} $
denotes the characteristic function of $\omega(t)$ (i.e.
$\chi_{\omega (t)}(x) = 1,$ when $x \in \omega(t)$ and
$\chi_{\omega(t)}(x) = 0,$ for $x \not\in \omega(t)$). In the
sequel, we will denote by $\tilde{y}_{y_0,u,r}$ the solution of
system $(S1)$ and by $y_{y_0,v,r}$ the solution of system $(S2)$. By the
linearity of $(S1)$ it is clear that $\tilde{y}_{y_0,u,r} =
\tilde{y}_{y_0} + \tilde{y}_{0,u,r}$, where $\tilde{y}_{y_0}$
denotes the unique solution of $(S1)$ with $u = 0$. If the subinterval
$\omega(t) = \omega$ for all $t$, (i.e. it is static), we will use
the notation $\tilde{y}_{y_0,u,\omega}$ instead of
$\tilde{y}_{y_0,u,r}$ (resp., $y_{y_0,v,\omega}$ instead of
$y_{y_0,v,r}$).

The control $u$ in $(S1)$ represents an external source acting
within $\omega(t)$. In turn, the multiplicative control $ v $ in
$(S2)$ can be viewed as controlling the reaction rate of the process
$(S2)$ within its support.

We will investigate how rich the sets of solutions to the above
controlled systems are at the fixed final time $T > 0$. To that end,
let us introduce the following notations:

\[  F_{1,y_0,mb}^+(T) = \{ \tilde{y}_{y_0,u,r}(x,T) : \ x \in (0,1), \
u \in L_+^2(Q_T) \ \mbox{ and } \ r \in PC[0,T]\}, \]

\begin{equation}  F_{2,y_0,mb}(T) = \{ y_{y_0,v,r}(x,T) :  \ x \in (0,1), \ v \in
L^\infty(Q_T) \ \mbox{ and } \ r \in PC[0,T] \}, \label{1.3}
\end{equation}
where $PC[0,T]$ denotes the set of functions $r:
[0,T]\longrightarrow [0, 1-l]$  such that $r$ is piecewise constant
with (at most) a finite number of discontinuities and $L_+^2 (Q_T) =
\{\phi \in L^2 (Q_T) : \phi(x,t) \geq 0 \;\; {\rm a.e. \;\; in } \;
Q_T\}.$

System $(S1)$ is said to be {\em approximately nonnegative
controllable at time $T > 0$} if for any given $y_0 \in L_+^2
(0,1)$, the set $F_{1,y_0,mb}^+(T)$ is dense in
$\tilde{y}_{y_0}(\cdot,T)+ L^2_+ (0,1)$. Taking into account the
linearity of the PDE, it is quite easy to verify
\begin{equation}
 F_{1,y_0,mb}^+(T) = \tilde{y}_{y_0}(\cdot,T) + F_{1,0,mb}^+(T). \label{1.2}
\end{equation}
Hence, the approximate nonnegative controllability at time $T > 0$
is clearly equivalent to the density of $F_{1,0,mb}^+(T)$ in $L^2_+
(0,1)$. Analogously, system $(S2)$ is said to be {\em approximately
nonnegative controllable at time $T > 0$} if for any given $y_0 \in
L_+^2 (0,1)$, $y_0 \not\equiv 0$, the set $ F_{2,y_0,mb}(T)$ is
dense in $ L^2_+ (0,1)$. For $y_0 \equiv 0$, by the uniqueness of
solution for the system $(S2)$, it is obvious that $F_{2,0,mb}(T) =
\{0\}$ and the aforementioned density property in $L^2_+ (0,1)$
clearly does not hold.

Formally, it is apparent that the PDE in $(S2)$ can be obtained from the
corresponding one in $(S1)$ by setting
\begin{equation} v(x,t) = \frac{u (x,t)}{y(x,t)}.
\label{1.35}
\end{equation}
This identity will be a key point to get our results. Of course,
these two problems are not equivalent in the sense that in $(S1)$ the
control $u$ is always nonnegative, while $v$ in $(S2)$ can be negative
as well, thus providing a ``qualitatively richer" set of solutions
for $(S2)$ compared to that of $(S1)$.

When the support of the control is the whole space domain, it was
shown in [8] that the set
\begin{equation}  F_{1,y_0}^+(T) = \{ \tilde{y}_{y_0,u,\omega}(x,T) : \  \ x \in (0,1), \ u \in
L_+^2(Q_T), \ \ \omega \equiv (0,1), \ \ \forall t \in (0,T) \} \label{1.4}
\end{equation}
is dense in $\tilde{y}_{y_0}(\cdot,T)+ L^2_+ (0,1)$. Furthermore, it
was shown in [9]-[11] that the corresponding property also holds for
the system $(S2)$: i.e., the set
\begin{equation}  F_{2,y_0}(T) = \{ y_{y_0,v,\omega}(x,T) :   \ x \in (0,1), \ v \in
L^\infty(Q_T), \ \ \omega \equiv (0,1), \ \ \forall t \in (0,T) \}
\label{1.5}
\end{equation}
is dense in $L^2_+(0,1)$. In fact, these results remain true for
more general linear and semilinear parabolic equations in several
space dimensions. In this paper, however, we are specifically
interested in the case when $\omega (t) $ is a proper subinterval of
$(0,1)$ at every moment of time.

Let us emphasize that the multiplicative controllability problem is
of nonlinear nature, even in the case of the linear equation like
$(S2)$, because the mapping $v \rightarrow y_{y_0,v,\omega}$ is highly
nonlinear: for instance, when $\omega = (0,1)$ and $v$ is constant,
it is well known that $y_{y_0,v,\omega}(x,t) = \tilde{y}_{y_0}(x,t) e^{-v t}$.
This means that we cannot use here the traditional duality argument,
which is the typical strategy to approach controllability problems
in the case of linear PDE with additive controls.

Among earlier works on the controllability of linear PDE by means of
bilinear controls we refer to the pioneering work \cite{3} by Ball,
Mardsen, and Slemrod. The paper \cite{13} further explored the
ideas of \cite{3}.

An extensive and thorough bibliography on controllability of
bilinear ODE is available. Let us just mention the survey
\cite{1}. Research in this area was seemingly originated in the 60s, on the
one hand, by the works of Kucera, who linked this area to the
Lie Algebra approach,  and, on the other hand, by the works of Rink and Mohler \cite{20} (see also \cite{12}), who pursued the qualitative
approach and numerous applications.

For recent works on global controllability of parabolic and
hyperbolic linear and semilinear PDE we refer to \cite{6}, \cite{8}-\cite{10} and
specially to the monograph \cite{11} and the references therein.

We also refer to the  works \cite{15}-\cite{16} (see also the references therein) on the issue of {\em
optimal} bilinear control for various PDE. A closely related issue
is {\em stabilization by means of bilinear controls}, see \cite{2} and \cite{18}.

\section{Main results}
\setcounter{equation}{0} First of all, we will show that {\em for
static} subdomain the above controllability property is out of the
question for both systems $(S1)$ and $(S2)$.

\begin{theorem}[{\bf Lack of controllability for static
support}] Let us assume that $y_0 \in L_+^2 (0,1)$, $y_0 \not\equiv
0$, and $\omega$ is a fixed proper subinterval of $(0,1)$. Then, for any $T > 0$, the
sets
\begin{equation}  F_{1,0,\omega}^+(T) = \{
\tilde{y}_{0,u,\omega}(x,T) :  \ x \in (0,1), \ u \in L_+^2(Q_T), \
\ \omega(t) \equiv \omega, \ \ \forall t \in (0,T)  \}, \label{2.2}
\end{equation}
and
\begin{equation}  F_{2,y_0,\omega}(T) = \{ y_{y_0,v,\omega}(x,T) :
\ x \in (0,1), \ v \in L^\infty(Q_T), \ \ \omega(t) \equiv \omega, \
\ \forall t \in (0,T) \} \label{2.1}
\end{equation}
are not dense in $L^2_+ (0,1)$. \label{T1}
\end{theorem}

Nevertheless, we can prove (in a constructive way) the following
``{\em local}'' controllability result.

\begin{theorem}[{\bf Local nonnegative controllability}]
Let us assume that $
\omega $ is a fixed proper subinterval of $(0,1)$. Then, for any $T > 0$, the set
\begin{equation}  F_{1,0,\omega}^+(T)_{|_{\omega}} = \{ \tilde{y}_{0,u,\omega}(x,T) :  \ x \in \omega, \ u \in
L_+^2(Q_T), \ \ \omega(t) \equiv \omega, \ \ \forall t \in (0,T) \},
\label{2.4}
\end{equation}
is dense in $L_+^2 (\omega)$. \label{T2}
\end{theorem}

If the control support is allowed to move with an additional
piecewise constant control $r(t)$, defining its position at time
$t$, we have the following {\em global controllability result}.

\begin{theorem}[{\bf Approximate controllability with mobile
support}] Let us assume that $y_0 \in L_+^2 (0,1)$, $y_0 \not\equiv
0$. Then, for any $T > 0$, the sets $F_{1,0,mb}^+(T)$ and $F_{2,y_0,mb}(T)$ are dense
in $L^2_+ (0,1)$. \label{T3}
\end{theorem}

In the proof of Theorem \ref{T1} we have faced with the following system:

$$ \left\{
\begin{array}{lcll}
y_{t}(x,t) &  =   & y_{xx}(x,t), & {\rm in} \; Q_T,
\\y(0,t) & = & u_0(t), &  {\rm in } \;\; (0,T),
\\ y(1,t) & = & u_1(t),   &  {\rm in } \;\; (0,T),
\\ y(x,0) &  = & y_0(x),  & {\rm in } \;\; (0,1),
\end{array}
\right. \eqno(S3)$$ with $u_0$ and $u_1$ acting as boundary
controls, meanwhile $y_0$ is fixed.

If we introduce the set
\begin{equation}  F_{y_0,b}(T) = \{ y_{y_0,u_0,u_1}(x,T) :  \ x \in (0,1), \ u_0, u_1 \in
L_+^\infty(0,T) \} \end{equation} where $y_{y_0,u_0,u_1}$ denotes
the unique solution of $(S3)$ defined by
transposition (see \cite[pp. 194-7]{17}), we have derived a related result on the lack of nonnegative
approximate boundary controllability with nonnegative boundary
controls, in the following terms:

\begin{theorem} Let us assume that $y_0 \in L_+^2 (0,1)$. Then, for any $T > 0$, the set $F_{y_0,b}(T)$ is not dense in $\tilde{y}_{y_0}(\cdot,T)+ L^2_+ (0,1)$.  \label{T4}
\end{theorem}

An essential tool for obtaining these results are the classical Maximum Principles for parabolic equations (both in the weak and strong form).
Since they will be cited very frequently along the paper and there are different versions depending on the framework, we recall them here:

\begin{theorem}[{\bf  Maximum Principles}]
Let us assume that $y_0 \in L^2(0,1)$, $u \in L^2(Q_T)$ and $v \in L^\infty(Q_T)$. Let $y \in L^2(0,T;H^1(0,1)) \cap C([0,T];L^2(0,1))$ be a solution of the problem
$$
\left\{
\begin{array}{lcll}
y_{t}(x,t) &  =   & y_{xx}(x,t) + v(x,t)y(x,t) + u(x,t), & {\rm
in} \; Q_T,
\\ y(x,0) &  = & y_0(x),  & {\rm in } \;\; (0,1).
\end{array}
\right.$$

\begin{itemize}
\item [i)] If $y_0 \geq 0$ in $(0,1)$, $u \geq 0$ in $Q_T$, $y(0,t)$ and $y(1,t) \geq 0$  in $(0,T)$, then $y(x,t) \geq 0$ a.e. in $Q_T$.
\item [ii)] If $u \equiv 0$, $v \leq 0$ in $Q_T$, $y(0,\cdot), y(1,\cdot) \in L^\infty(0,T)$ and $y_0 \in L^\infty (0,1)$, then $y(x,t) \leq K$ a.e. in $Q_T$, with
\[ K = \max\{ \|(y(0,\cdot))^+\|_{L^\infty(0,T)}, \|(y(1,\cdot))^+\|_{L^\infty(0,T)}, \|(y_0)^+\|_{L^\infty(0,1)} \},\]
where $z^+ = max\{z,0\}.$
\item [iii)] ({\bf Strong}) If $y \in L^2(0,T;H_0^1(0,1))$, $y_0 \geq 0$ in $(0,1)$ ($y_0 \not\equiv 0$), $u \geq 0$ in $Q_T$, $v = v(x) \leq 0$ in $(0,1)$, then $y(x,t) >  0$ in $(0,1)$ for all $t \in (0,T]$.

\end{itemize}
\label{MP}
\end{theorem}

For the proofs, we refer (among others) to \cite[Chapter III, section 7, pp. 181-191]{14} and \cite[pp. 21-23]{19}. See also \cite[Chapter I, section 2, pp. 11-25]{14} for related results in the classical framework.

\section{Proof of Theorems \ref{T1}, \ref{T2} and \ref{T4}.}
\setcounter{equation}{0} We begin with the proof of Theorem \ref{T1} for
system $(S1)$. Our goal is to construct a simple explicit
counterexample.

\subsection{Proof of Theorem \ref{T1}  for system $(S1)$.}

Since $u \in L^2_+(Q_T)$, it is well known (see for instance \cite[Theorem 4.2, p. 160]{14}), that each solution $\tilde{y}_{0,u,\omega}$ belongs to
$L^2(0,T;H_0^1(0,1)) \cap C([0,T];L^2(0,1))$. Moreover, due to the Theorem \ref{MP}-i)
\[ \tilde{y}_{0,u,\omega}(x,t)
\geq 0 \ \ \mbox{a.e. } x \in (0,1), \ \forall \ t \in [0,T].\]
Consequently, $F_{1,0,\omega}^+(T) \subset L^2_{+}(0,1)$. We will
show that $F_{1,0,\omega}^+(T)$ is not dense in $L^2_{+}(0,1)$.

Since $\omega $ is a fixed proper subinterval of $(0,1)$, there will
exist a sufficiently large natural number $m \geq 2$ such that
$\omega \subset (\frac{1}{m},1)$. Let us fix that $m$ for the rest
of the proof. It is now straightforward to check that if
$$ p(x,t) =
\varphi(x)\exp{\left(\frac{m^2\pi^2(t-T)}{(m-1)^2}\right)},$$ with
$$ \varphi(x) = \left\{\begin{array}{ll}\sin{(m\pi x)}, & x \in [0,\frac{1}{m}), \\
(1-m)\sin{(\frac{\pi(mx-1)}{m-1})}, & x \in
[\frac{1}{m},1],\end{array}\right.
$$ $p \in C^2(\overline{Q_T})$ and it is the unique classical solution of the adjoint problem
 \begin{equation}
\left\{
\begin{array}{lcll}
-p_{t}(x,t) &  =   & p_{xx}(x,t) + h(x,t),& {\rm in} \; Q_T,
\\p(0,t) & = & p(1,t)  = 0,   &  {\rm in } \;\; (0,T),
\\ p(x,T) &  = & \varphi(x),  & {\rm in } \;\; (0,1),
\end{array}
\right. \label{3.1}
\end{equation}
where
$$ h(x,t) = \left\{\begin{array}{ll}\frac{(m-2)m^3\pi^2}{(m-1)^2}\sin{(m\pi x)}\exp{\left(\frac{m^2\pi^2(t-T)}{(m-1)^2}\right)}, & x \in [0,\frac{1}{m}), \\
0, & x \in [\frac{1}{m},1].\end{array}\right.$$ Furthermore,

\begin{equation}h(x,t) \geq 0 \;\;\;\; {\rm in} \;\; \overline{Q_T},
\label{3.1a}
\end{equation}

\begin{equation}p(x,t) \leq 0 \;\;\;\; {\rm in} \;\; \left[\frac{1}{m},1\right] \times [0,T].
\label{3.2}
\end{equation}

Let us now pick up the state in $L^2_{+}(0,1)$
\[ y_d(x) = \max{\{\varphi(x),0\}}. \]

We will see that $y_d$ is ``unachievable''. Arguing by
contradiction, let us suppose that there exists a sequence
$\{u_k\}_{k=1}^\infty \subset L^2_{+}(Q_T)$ such that:
\[ \tilde{y}_{0,u_k,\omega}(\cdot, T) \longrightarrow y_d \ \ \mbox{ in } L^2(0,1) \;\; {\rm as} \;  k \rightarrow + \infty. \]

Multiplying $(S1)$ (with $u = u_k$, $\omega(t) = \omega$ and $y_0 =
0$) by $p$, further integration by parts, combined with
(\ref{3.1})-(\ref{3.2}) provides:

\[ \int_0^1 \tilde{y}_{0,u_k,\omega}(x,T)\varphi(x) dx \leq \int_0^1 \tilde{y}_{0,u_k,\omega}(x,T)\varphi(x) dx + \int_{Q_T} h(x,t) \tilde{y}_{0,u_k,\omega}(x,t) dx dt = \]
\begin{equation}
= \int_{\omega \times (0,T)} u_k(x,t) p(x,t) dx dt \leq 0 \; \;\;\;
\forall k. \label{3.3}
\end{equation}

Thus, we arrive at the contradiction
\[ 0 < \int_0^{\frac{1}{m}} (\sin{(m\pi x)})^2
dx = \int_0^1 y_d(x) \varphi(x) dx = \lim_{k} \int_0^1
\tilde{y}_{0,u_k,\omega}(x,T) \varphi(x) dx \leq 0.\] This ends the
proof of Theorem \ref{T1} for system $(S1)$.\rule{2mm}{2mm}

\subsection{Proof of Theorem \ref{T2}.}
It is enough to show that solutions to system $(S1)$ (with $y_0 = 0$
and fixed $\omega$) at time $T$ can approximate in $L_+^2(\omega)$
any element of $H_0^1 (\omega) \bigcap L_+^2 (\omega)$. This can be
shown by a simple modification of the proof given in \cite{5} for the
case $\omega = (0,1)$. Nonetheless, the argument used there does not
provide any insight on what controls are used to achieve the
desirable steering. Therefore, here we give a different constructive
proof of Theorem \ref{T2}, through Fourier series expansions, which
can be used to prove the aforementioned result in \cite{5} as well and will also serve us as a
tool to prove the Theorem \ref{T3} below.

Consider any element  $y_d \in H_0^1 (\omega) \bigcap L_+^2
(\omega)$. Extend it by zero to the whole interval $(0,1)$ and, for simplicity of
notations, denote the resulting function by $ y_d$ again. Clearly,
it lies in $H_0^1(0,1) \bigcap L_+^2 (0,1)$. Given any $\delta \in (0, T)$, we define
\begin{equation} u_\delta(x,t) = \left\{
\begin{array}{ll}
0, & \ t \in [0,T - \delta),
\\ \frac{1}{\delta}  y_d (x),    & \ t \in [T - \delta, T].
\end{array}
\right. \label{3.4} \end{equation}
Taking now advantage of the Fourier series expansion, we know
\begin{equation} y_d(x) = 2\sum_{k = 1}^\infty a_k \sin (k \pi x), \ \ \mbox{ in } \ H_0^1(0,1), \label{3.4-1} \end{equation}
with $a_k = \int_0^1 y_d (s) \sin (k \pi s) ds$.

The classical method of
separation of variables provides us the following expression for the
corresponding solution of $(S1)$ with $y_0 = 0$
\begin{equation}
\tilde{y}_{0,u_\delta,\omega}(x,t) = \left\{
\begin{array}{ll}
0, & \ t \in [0,T - \delta],
\\ \\ \frac{2}{\delta}  \sum_{k
= 1}^\infty \frac{1- e^{-\pi^2 k^2 (t-T+\delta)}}{\pi^2 k^2} a_k \sin (k \pi
x),    & \ t \in [T - \delta, T].
\end{array}
\right. \label{3.5}\end{equation} Hence, we derive that
\begin{equation} \parallel \tilde{y}_{0,u_\delta,\omega} (\cdot, T) - y_d \parallel^2_{H_0^1 (0,1)} \; =
2\pi^2\sum_{k = 1}^\infty a_k^2  k^2 \left(\frac{1-e^{-\pi^2
k^2 \delta}}{\pi^2 k^2 \delta}-1\right)^2, \label{3.6} \end{equation} where we endowed the Sobolev space $H_0^1 (0,1)$ with the norm
$$
\parallel \phi \parallel_{H_0^1 (0,1)} \; = \; \left( \int_0^1 \phi_x(x)^2 dx \right)^{1/2}.
$$

Let us show that the right-hand side of (\ref{3.6}) tends to zero
as $ \delta \rightarrow 0^+$. To that end, we introduce
the auxiliary function \[ \psi(r) =
\left(\frac{1-e^{-r}}{r}-1\right)^2, \ \ r \geq 0.
\] It can be verified in a straightforward way that $\psi$ is a $C^\infty$ function, strictly increasing in $(0,+\infty)$, $\psi(0)=0$ and $\psi([0,+\infty)) \subset [0,1)$.

Taking into account that $y_d \in H_0^1 (0,1)$, given any $\epsilon
>0$, we can find a natural number $N$ such that,
\begin{equation}
\sum_{k = N+1}^\infty a_k^2 k^2 \; \leq \;
\frac{\epsilon}{4\pi^2}. \label{3.7}
\end{equation} Using now that $\psi$ is continuous at $0$ and $\psi(0) = 0$, we can select $
\hat{\delta} >0 $ small enough to guarantee that
\begin{equation}
\psi(\pi^2 N^2 \delta)\sum_{k = 1}^N a_k^2 k^2 \leq \;
\frac{\epsilon}{4\pi^2}, \ \ \forall \ \delta \in [0,\hat{\delta}].\label{3.8}
\end{equation}

Combining previous estimates with the fact that $\psi(\pi^2 k^2 \delta) < 1$ for all $k$, we get
that
\begin{equation}
\parallel \tilde{y}_{0,u_\delta,\omega} (\cdot, T) - y_d \parallel^2_{H_0^1 (0,1)} \; \leq \;
2\pi^2 \left(\psi(\pi^2 N^2 \delta) \sum_{k = 1}^N a_k^2 k^2 + \sum_{k =
N+1}^\infty a_k^2 k^2\right) \leq \epsilon, \label{3.9}
\end{equation}
for all $\delta \in (0,\hat{\delta}]$. This ends the proof of Theorem \ref{T2}, but for proving later Theorem \ref{T3} corresponding to
system $(S1)$, it is convenient to make a small modification to
previous construction as follows: instead of the control $u_\delta(x,t)$
given by (\ref{3.4}), we can take a
control with value zero at the end of the process. More precisely, we can consider for instance
\begin{equation} \hat{u}_\delta(x,t) = \left\{
\begin{array}{ll}
0, & \ t \in (0,T - 2\delta),
\\\frac{1}{\delta}  y_d (x),    & \ t \in [T - 2\delta, T - \delta),
\\ 0, & \ t \in [T - \delta,T].
\end{array}
\right. \label{3.11r}
\end{equation}

Applying the method of separation of variables in each time interval, it is not difficult to see here that the solution of $(S1)$ with $y_0 = 0$ is given by

\begin{equation}
\tilde{y}_{0,\hat{u}_\delta,\omega}(x,t) = \left\{
\begin{array}{ll}
0, & \ t \in [0,T - 2\delta],
\\ \\ \frac{2}{\delta}  \sum_{k
= 1}^\infty \frac{1- e^{-\pi^2 k^2 (t-T+2\delta)}}{\pi^2 k^2} a_k \sin (k \pi
x),    & \ t \in [T - 2\delta, T -\delta],
\\ \\ \frac{2}{\delta}  \sum_{k
= 1}^\infty \frac{1- e^{-\pi^2 k^2 \delta}}{\pi^2 k^2}e^{-\pi^2 k^2 (t-T+\delta)} a_k \sin (k \pi
x),    & \ t \in [T - \delta, T],
\end{array}
\right. \label{3.12a}\end{equation} and therefore, arguing as before, we arrive to
\begin{equation} \parallel \tilde{y}_{0,\hat{u}_\delta,\omega} (\cdot, T) - y_d \parallel^2_{H_0^1 (0,1)} \; =
2\pi^2\sum_{k = 1}^\infty a_k^2 k^2 \left(\frac{e^{-\pi^2
k^2 \delta}-e^{-2\pi^2
k^2 \delta}}{\pi^2 k^2 \delta}-1\right)^2, \label{3.12b} \end{equation}

We can finish as above, seeing that the right-hand side of (\ref{3.12b}) tends to zero
as $ \delta \rightarrow 0^+$ with the help of the new auxiliary function \[ \hat{\psi}(r) =
\left(\frac{e^{-r}-e^{-2r}}{r}-1\right)^2, \ \ r \geq 0,
\] that satisfies the same properties than $\psi$.\rule{2mm}{2mm}

\subsection{Proof of Theorem \ref{T1} for system $(S2)$.}
Let $(\alpha,\beta)$ be any subinterval of $(0,1) \backslash \omega$
and let us pick out any $y_d \in L^2_+(0,1)$ such that $y_d \equiv 0$ in $(\alpha, \beta)$.
Arguing by contradiction, let us assume that there exist
$\{v_k\}_{k=1}^\infty \subset L^\infty(Q_T)$ such that
\[ y_{y_0,v_k,\omega}(\cdot,T) \rightarrow y_d  \ \ \mbox{ in } L_+^2(0,1) \;\; {\rm as} \;  k \rightarrow + \infty. \]
In particular, we have
\[ y_{y_0,v_k,\omega}(\cdot,T)_{|_{(\alpha,\beta)}} \rightarrow 0  \ \ \mbox{ in } L_+^2(\alpha,\beta) \;\; {\rm as} \;  k \rightarrow + \infty. \]

The restriction of each $y_{y_0,v_k,\omega}$ to the interval $(\alpha,\beta)$ can be viewed as the solution
of one system of the type
$$ \left\{
\begin{array}{lcll}
y_{t}(x,t) &  =   & y_{xx}(x,t), & {\rm in} \; (\alpha,\beta) \times (0,T),
\\y(\alpha,t) & = & u_0(t), &  {\rm in } \;\; (0,T),
\\ y(\beta,t) & = & u_1(t),   &  {\rm in } \;\; (0,T),
\\ y(x,0) &  = & y_0(x),  & {\rm in } \;\; (\alpha,\beta),
\end{array}
\right. \eqno(S3a)$$
with nonnegative boundary values at the extremes $\alpha$ and $\beta$, because
$y_{y_0,v_k,\omega}$ is nonnegative thanks to the Theorem \ref{MP}-i).

On the other hand, if we denote by $\hat{y}$ the unique solution of the system
$$\left\{
\begin{array}{lcll}
y_{t}(x,t) &  =   & y_{xx}(x,t), & {\rm in} \; (\alpha, \beta) \times (0,T),
\\y(\alpha,t) & = & y(\beta,t) = 0, &  {\rm in } \;\; (0,T),
\\ y(x,0) &  = & y_0(x),  & {\rm in } \;\; (\alpha, \beta),
\end{array}
\right. \eqno(S3b)$$
as a consequence of the Theorem \ref{MP}-i) we have that
\[ y_{y_0,v_k,\omega}(x,t)_{|_{(\alpha,\beta)}}- \hat{y}(x,t) \geq 0, \ \ \ a. e. \ (x,t) \in (\alpha, \beta) \times (0,T),\]
and thanks to Theorem \ref{MP}-iii)
\[ \hat{y}(x,T) > 0, \ \ \ a. e. \ x \in (\alpha, \beta).\]
The contradiction arrives combining all the previous conditions in the form
\[ 0 = \lim_{k \rightarrow + \infty} \|y_{y_0,v_k,\omega}(\cdot,T)\|_{L^2(\alpha,\beta)} \geq \|\hat{y}(\cdot,T)\|_{L^2(\alpha,\beta)} > 0. \rule{2mm}{2mm}\]

\subsection{Proof of Theorem \ref{T4}.}

In the previous proof, we have faced with a problem of the type $(S3)$ (more precisely, $(S3a)$).
Now, we will show that the nonnegative controllability property for these
problems is out of the question, even if the obstruction phenomenon due to the initial datum is removed.

By the linearity of the problem, it is well known that $y_{y_0,u_0,u_1} = \tilde{y}_{y_0} + y_{0,u_0,u_1},$
and the following relation holds
\begin{equation}
F_{y_0,b}(T) = \tilde{y}_{y_0}(\cdot,T)+ F_{0,b}(T). \label{2.6}
\end{equation}
Therefore, the conclusion of Theorem  \ref{T4} is equivalent to say that $F_{0,b}(T)$ is not dense in $L^2_+
(0,1)$ and this is what we will prove.

At this point, let us remind that (in general) the solution $y_{0,u_0,u_1}$ must be defined by the
transposition method (see \cite[pp. 194-7]{17}). In our situation, its expression is given by
\begin{equation}
y_{0,u_0,u_1}(x,t) = \sum_{k=1}^\infty b_k(t) \sqrt{2}\sin{(k\pi x)},
\label{3.13-1}
\end{equation}
where
\begin{equation}
b_k(t) = \sqrt{2}k \pi \int_0^t e^{-k^2\pi^2(t-s)}(u_0(s)+(-1)^{k+1}u_1(s))ds.
\label{3.13-2}
\end{equation}
It is not obvious at all that $y_{0,u_0,u_1}(\cdot,T) \in L^2(0,1)$. In fact, in \cite[p. 202]{17}, it is constructed an example where $y_{0,u_0,u_1}(\cdot,T) \not\in L^2(0,1)$, by taking $u_0(t) = 0$ and $u_1(t) = \frac{1}{\sqrt[4]{T-t}}$.
The difficulty here is that $u_1 \in L^2(0,T)$, but not in $L^\infty(0,T)$. In our case, using (\ref{3.13-1})-(\ref{3.13-2}), we can estimate
\begin{equation}
\int_0^1|y_{0,u_0,u_1}(x,T)|^2 dx = \sum_{k=1}^\infty |b_k(T)|^2 \leq C_1\sum_{k=1}^\infty \frac{1}{k^2} < \infty,
\label{3.13-3}
\end{equation}
with $C_1 = \frac{2}{\pi^2}\left(\|u_0\|_{L^\infty(0,T)}+\|u_1\|_{L^\infty(0,T)}\right)^2$.

Moreover, applying the Theorem \ref{MP}-i) we derive that $y_{0,u_0,u_1}(x,T) \geq 0,$ a. e. $x \in (0,1).$ Consequently, $F_{0,b}(T)$ is included in $L^2_+(0,1)$. Let us show that it is not dense in $L^2_+(0,1)$, arguing by contradiction.
To this end, it can be checked by direct calculation that
\[ p(x,t) = - \exp{(9 \pi^2(t-T))}\sin{(3
\pi x)}\] is the unique classical solution of the adjoint problem
\begin{equation}
\left\{
\begin{array}{lcll}
-p_{t}(x,t) &  =   & p_{xx}(x,t), & {\rm in} \; Q_T,
\\p(0,t) & = & p(1,t)  = 0,   &  {\rm in } \;\; (0,T),
\\ p(x,T) &  = & - \sin{(3\pi x)},  & {\rm in } \;\; (0,1).
\end{array}
\right. \label{3.13}
\end{equation}
Moreover, it satisfies $p_x(0,t) < 0, \ p_x(1,t) > 0$ in $(0,T)$. Now, we select
\[ y_d(x) = \max{\{p(x,T),0\}}, \;\; x \in (0,1), \]
as the target state in $L^2_{+}(0,1)$.

If $F_{0,b}(T)$ is dense in $L^2_{+}(0,1)$, there exist two
sequences $\{u_{0k}\}_{k=1}^\infty, \{u_{1k}\}_{k=1}^\infty$ in
$L^\infty_{+}(0,T)$ such that for the corresponding solutions
$\{y_k\}_{k = 1}^\infty$ of $(S3)$ with $y_0 = 0$ (namely, $y_k =
y_{0,u_{0k},u_{1k}}$), we have:

\[ y_k(\cdot, T) \longrightarrow y_d \ \ \mbox{ in } L_+^2(0,1) \;\; {\rm as} \;  k \rightarrow + \infty. \]

Multiplying by $p$ the PDE satisfied by $y_k$ and integrating by
parts, we get:
\[ \int_0^1 y_k(x,T)p(x,T) dx = \int_0^T u_{0k}(t)p_x(0,t) dt -
\int_0^T u_{1k}(t)p_x(1,t) dt \leq 0. \]
(Alternatively, previous equality can be derived from the expressions (\ref{3.13-1})-(\ref{3.13-2}) for each $y_k$.)
Thus, we arrive to the contradiction
\[ \int_{1/3}^{2/3} (\sin{(3\pi x)})^2
dx = \int_0^1 y_d(x) p(x,T) dx = \lim_{k \rightarrow \infty}
\int_0^1 y_k(x,T) p(x,T) dx \leq 0. \rule{2mm}{2mm} \]

\section{Proof of Theorem \ref{T3}.}
\setcounter{equation}{0}

\subsection{Proof of Theorem \ref{T3} for system $(S1)$.}
To show that $F_{1,0,mb}^+(T)$ is dense in $L^2_+ (0,1)$, we will
use the technique developed in the proof of Theorem \ref{T2}. First,
let us fix $M$ the smallest natural number satisfying $M \cdot l
\geq 1$, where $l$ is the (fixed) length of the mobile subinterval
$\omega(t)$ (see (\ref{1.1})).

Given any element $y_d \in L^2_+(0,1)$, we can decompose it into $M$
``pieces" localized in disjoint subintervals of length (at most)
$l$, as follows:
\[y_{d_1}(x) = \left\{
\begin{array}{ll}
y_d(x), & \ x \in (0,l),
\\ 0,  & x \in (l,1),
\end{array}
\right.
\]
\[y_{d_j}(x) = \left\{
\begin{array}{ll}
y_d(x), & \ x \in ((j-1)l,lj),
\\ 0,  &  \mbox{otherwise},
\end{array}
\right. \] for $j =2,\ldots,M-1$,

\begin{equation} y_{d_M}(x) = \left\{
\begin{array}{ll}
y_d(x), & \ x \in ((M-1)l,1),
\\ 0,  & \mbox{otherwise}.
\end{array}
\right. \label{4.1}
\end{equation}
 Clearly, each $y_{d j} \in
L^2_+(0,1)$ and
\[ y_d(x) = \sum_{j=1}^M y_{d_j}(x), \ \ \ \mbox{a.e. } \ x \in
(0,1).\]

Moreover, given $\epsilon > 0$, there exist $\hat{y}_{d_j} \in
H_0^1(0,1) \cap L^2_+(0,1)$ for $j =1,\ldots,M$, with its support $\omega_j$ contained in the corresponding
one of each $y_{d j}$ and verifying

\begin{equation}
\sum_{j=1}^M \parallel \hat{y}_{d_j} -y_{d_j}
\parallel_{L^2(0,1)} \leq \frac{\epsilon}{2}.  \label{4.2}
\end{equation}

Now, we can argue ``by parts", using the constructions (\ref{3.11r}) and (\ref{3.4}) established in the proof of Theorem \ref{T2}, as follows: given $\hat{y}_{d_1}$, there exists $\delta_1 \in (0,\frac{T}{2})$ such that if we take the control

\begin{equation} \hat{u}_1(x,t) = \left\{
\begin{array}{ll}
0, & \ t \in (0,T - 2\delta_1),
\\\frac{1}{\delta_1}  \hat{y}_{d_1} (x),    & \ t \in (T - 2\delta_1, T - \delta_1),
\\ 0, & \ t \in (T - \delta_1,T),
\end{array}
\right. \label{4.3}
\end{equation}
it is satisfied
\begin{equation}
\parallel \tilde{y}_{0,\hat{u}_1,\omega_1} (\cdot, T) - \hat{y}_{d_1} \parallel_{L^2(0,1)} \; \leq  \; \frac{\epsilon}{2M}. \label{4.4}
\end{equation}
Arguing now with $\hat{y}_{d_2}$, there exists $\delta_2 \in
\left(0,\frac{\delta_1}{2}\right)$ such that if we take the control

\begin{equation} \hat{u}_2(x,t) = \left\{
\begin{array}{ll}
0, & \ t \in (0,T - 2\delta_2),
\\\frac{1}{\delta_2}  \hat{y}_{d_2} (x),    & \ t \in (T - 2\delta_2, T - \delta_2),
\\ 0, & \ t \in (T - \delta_2,T),
\end{array}
\right. \label{4.5}
\end{equation}
it is satisfied
\begin{equation}
\parallel \tilde{y}_{0,\hat{u}_2,\omega_2} (\cdot, T) - \hat{y}_{d_2} \parallel_{L^2(0,1)} \; \leq  \; \frac{\epsilon}{2M}. \label{4.6}
\end{equation}

Iterating this process, we obtain $\delta_j \in \left(0,\frac{\delta_{j-1}}{2}\right),$ $j = 1,\ldots,M$, and controls

\begin{equation} \hat{u}_j(x,t) = \left\{
\begin{array}{ll}
0, & \ t \in (0,T - 2\delta_j),
\\\frac{1}{\delta_j}  \hat{y}_{d_j} (x),    & \ t \in (T - 2\delta_j, T - \delta_j),
\\ 0, & \ t \in (T - \delta_j,T),
\end{array}
\right. \label{4.7}
\end{equation}
$j = 1,\ldots,M-1$, and

\begin{equation} \hat{u}_M(x,t) = \left\{
\begin{array}{ll}
0, & \ t \in (0,T - \delta_M),
\\\frac{1}{\delta_M}  \hat{y}_{d_M} (x),    & \ t \in (T - \delta_M, T),
\end{array}
\right.
\end{equation}
for which
\begin{equation}
\parallel \tilde{y}_{0,\hat{u}_j,\omega_j} (\cdot, T) - \hat{y}_{d_j} \parallel_{L^2(0,1)} \; \leq  \; \frac{\epsilon}{2M}, \ \ \ j = 1,\ldots,M. \label{4.8}
\end{equation}

We finish the proof by taking
\begin{equation}
u(x,t) = \sum_{j=1}^M \hat{u}_j(x,t), \ \ \ \ \ \tilde{y}(x,t) = \sum_{j=1}^M \tilde{y}_{0,\hat{u}_j,\omega_j}(x,t)\label{4.9}
\end{equation}
and noticing that by the linearity
\[ \parallel \tilde{y}(\cdot, T) - y_d \parallel_{L^2(0,1)}
\leq \; \parallel \tilde{y}(\cdot, T) - \sum_{j=1}^M
\hat{y}_{d_j} \parallel_{L^2(0,1)} + \parallel \sum_{j=1}^M
\hat{y}_{d_j}-y_d
\parallel_{L^2(0,1)} \leq \]
\begin{equation}
\leq \sum_{j=1}^M \parallel \tilde{y}_{0,\hat{u}_j,\omega_j} (\cdot, T) -
\hat{y}_{d_j}
\parallel_{L^2(0,1)} + \sum_{j=1}^M \parallel \hat{y}_{d_j} -y_{d_j}
\parallel_{L^2(0,1)} \leq  \; \epsilon. \label{4.10}
\end{equation}

Let us remark that in each time subinterval $(T - \delta_j, T - \delta_{j+1})$,
the value of the control $u(x,t)$ coincides with the value of $\hat{u}_{j+1}(x,t)$ and hence its support is included in $\omega(t) = (r(t),r(t)+l)$, with $r(t) \in PC[0,T]$ given by

\begin{equation} r(t) = \left\{
\begin{array}{ll}
0, & \ t \in (0,T - \delta_1),
\\ l,    & \ t \in (T - \delta_1, T - \delta_2),
\\ 2l,    & \ t \in (T - \delta_2, T - \delta_3),
\\ & \vdots
\\ jl,    & \ t \in (T - \delta_j, T - \delta_{j+1}),
\\ & \vdots
\\ (M-2)l,    & \ t \in (T - \delta_{M-2}, T - \delta_{M-1}),
\\ 1-l,    & \ t \in (T - \delta_{M-1}, T).
\end{array}
\right. \label{4.11}
\end{equation}

This implies that $u(x,t) = u(x,t)\chi_{\omega(t)}(x)$ and
$\tilde{y} = \tilde{y}_{0,u,r}$, as desired.\rule{2mm}{2mm}

\subsection{Proof of Theorem \ref{T3} for system $(S2)$.}

We will prove the result in two steps. First, following the
notations of Section 1 (compare with (\ref{1.3})), we introduce the
set
\begin{equation}  F_{2,y_0,mb}^+(T) = \{ y_{y_0,v,r}(x,T) :  \ x \in (0,1), \ v \in
L^\infty_+(Q_T) \ \mbox{ and } \ r \in PC[0,T] \}. \label{5.1}
\end{equation}

Applying twice the Theorem \ref{MP}-i) (first to $y_{y_0,v,r}$ and later to the difference $y_{y_0,v,r}-\tilde{y}_{y_0}$, treating the bilinear term as a free one), it is well known that
\[ y_{y_0,v,r}(x,T) \geq \tilde{y}_{y_0}(x,T), \ \ \mbox{a.e.} \ x \in (0,1), \ \forall v \in
L^\infty_+(Q_T),  \ \forall r \in PC[0,T]. \]

\subsubsection{Step 1.} We will prove that the system $(S2)$ is
approximately nonnegative controllable at time $T$, in the following sense:

\begin{theorem}[Nonnegative multiplicative mobile control] Let us assume that $y_0 \in L_+^2 (0,1)$, $y_0 \not\equiv 0$. Then,
for any $T > 0$, the set $F_{2,y_0,mb}^+(T)$ is dense in $\tilde{y}_{y_0}(\cdot,T) + L^2_+ (0,1)$. \label{T5}
\end{theorem}
{\bf Proof of Theorem \ref{T5}.} Pick out a generic element $y_d \in L^2_+ (0,1)$.
As a consequence of Theorem \ref{T3} for system $(S1)$, given
$\epsilon > 0$, there exist $u \in L^2_+(Q_T)$ and $r \in PC[0,T]$ such that
\begin{equation}
\parallel \tilde{y}_{0,u,r} (\cdot, T) - y_d
\parallel_{L^2(0,1)} < \epsilon. \label{5.2}
\end{equation}

By well known density results, we can suppose, without loss of generality, that $u \in
C_+(\overline{Q_T})$ (that is, $u$ is continuous and nonnegative in $\overline{Q_T}$) with compact support $K$ included
in $(0,1)\times(\delta,T) \subset Q_T,$ for some $\delta >0$.

Using classical regularity results (see \cite[Theorem 10.2, p. 140]{7} and \cite[Theorem 9.1, p. 341-2]{14} it is
also known that $\tilde{y}_{y_0} \in C^\infty(Q_T) \cap C_+([0,1]\times [\delta,T])$.

Now, let us consider $y = \tilde{y}_{y_0,u,r}$ and

\begin{equation} v(x,t) = \left\{
\begin{array}{cl}
\frac{u(x,t)}{y(x,t)},    & \ (x,t) \in K,
\\ 0, & \ (x,t) \in Q_T \setminus K.
\end{array}
\right. \label{5.3}
\end{equation}

By the Theorem \ref{MP}-i), it is clear that $y \geq 0,$
a.e. in $Q_T$. We will prove below that $v \in L^\infty_+(Q_T)$. Using
this fact, it is straightforward to see that $y = y_{y_0,v,r}$ (i.e. $y$
is the solution of system $(S2)$ corresponding to $y_0,v,r$), because
$u(x,t) \chi_{\omega(t)}(x) = v(x,t)y(x,t)\chi_{\omega(t)}(x)$ in
$Q_T$, and

\[ \parallel y_{y_0,v,r} (\cdot, T) - (\tilde{y}_{y_0}(\cdot,T)+y_d)
\parallel_{L^2(0,1)}= \]
\begin{equation}
= \parallel \tilde{y}_{y_0,u,r} (\cdot, T) - (\tilde{y}_{y_0}(\cdot,T)+y_d)
\parallel_{L^2(0,1)}= \parallel \tilde{y}_{0,u,r} (\cdot, T) - y_d
\parallel_{L^2(0,1)} < \epsilon, \label{5.4}
\end{equation}
as we were looking for.

Hence, it remains to show that $v \in L^\infty_+(Q_T)$. Since
$\tilde{y}_{y_0} \in C_+([0,1]\times [\delta,T])$,
$\tilde{y}_{y_0}\geq 0$ in $Q_T$ and $K$ is a compact set, there is
a constant $\rho \geq 0$ such that $\tilde{y}_{y_0} \geq \rho$ in
$K$. In fact, we can guarantee that $\rho >0$: if there exists
$(x_1,t_1) \in K \subset (0,1)\times (\delta,T)$ such that
$\tilde{y}_{y_0}(x_1,t_1) = \rho = 0$, thanks to the Theorem \ref{MP}-iii), we conclude that $\tilde{y}_{y_0}
\equiv 0$ in $[0,1] \times [0,t_1]$ and, consequently, $y_0 \equiv 0$ in $(0,1)$,
which it is excluded by hypothesis. Then, we get
\[ \parallel v \parallel_{L^\infty(Q_T)} \leq \frac{\parallel u
\parallel_{L^\infty(K)}}{\rho}. \rule{2mm}{2mm}\]

\subsubsection{Step 2.} We want to approximate any element $y_d \in
L^2_+ (0,1)$, even if $y_d \not\in \tilde{y}_{y_0}(\cdot,T) + L^2_+
(0,1)$. Roughly speaking, the main idea here is to take very large negative
controls at the beginning of the process and move the control support across
the space domain, in order to reduce the lower barrier that
represents $\tilde{y}_{y_0}(\cdot,T)$, making it small enough.

For that purpose, we fix $M$ (as previously) the
smallest natural number satisfying $M \cdot l \geq 1$, where $l$ is
the length of the mobile subinterval $\omega(t)$. Given some intermediate times $0 \leq T_1 \leq T_2 \leq \ldots \leq T_M < T$, let us
denote by $y_1(x,t)$ the unique solution of the following system
\begin{equation} \left\{
\begin{array}{lcll}
y_{t}(x,t) &  =   & y_{xx}(x,t) - m_1 y(x,t) \chi_{(0,l)}(x), & {\rm
in} \; (0,1)\times(0,T_1),
\\y(0,t) & = & y(1,t) =  0,   &  {\rm in } \;\; (0,T_1),
\\ y(x,0) &  = & y_0 (x),  & {\rm in } \;\; (0,1),
\end{array}
\right. \label{5.5}
\end{equation}
where $m_1 \geq 0$ will be carefully chosen later. Assuming now that $y_{j-1}$ is already known,
$y_j$ will denote the solution of the system
\begin{equation} \left\{
\begin{array}{lcll}
y_{t}(x,t) &  =   & y_{xx}(x,t) - m_j y(x,t) \chi_{((j-1)l,jl)}(x),
& {\rm in} \; (0,1)\times(T_{j-1},T_j),
\\y(0,t) & = & y(1,t) =  0,   &  {\rm in } \;\; (T_{j-1},T_j),
\\ y(x,T_{j-1}) &  = & y_{j-1}(x,T_{j-1}),  & {\rm in } \;\; (0,1),
\end{array}
\right. \label{5.6}
\end{equation}
with $m_j \geq 0$ for $j = 1,\ldots,M-1$ also to be chosen. Finally, $y_M$ will denote the solution of the problem

\begin{equation} \left\{
\begin{array}{lcll}
y_{t}(x,t) &  =   & y_{xx}(x,t) - m_M y(x,t) \chi_{((M-1)l,1)}(x),
& {\rm in} \; (0,1)\times(T_{M-1},T_M),
\\y(0,t) & = & y(1,t) =  0,   &  {\rm in } \;\; (T_{M-1},T_M),
\\ y(x,T_{M-1}) &  = & y_{M-1}(x,T_{M-1}),  & {\rm in } \;\; (0,1),
\end{array}
\right. \label{5.6-1}
\end{equation}
with $m_M \geq 0$ to be chosen. The following result will be essential in this section:

\begin{theorem} With previous notations, given $\epsilon >0$, there exist nonnegative constants $m_1,\ldots, m_M$
 and intermediate times $0 \leq T_1 \leq T_2 \leq \ldots \leq T_M < T$ such that
\[\|y_M(\cdot,T_M)\|_{L^2(0,1)} \leq
\frac{\epsilon}{2}.\] \label{T6}
\end{theorem}
We postpone its proof to the next sections, but we will conclude now the proof of Theorem \ref{T3}. Given $\hat{y}_1(x,t)$ the unique solution of the system
\begin{equation} \left\{
\begin{array}{lcll}
y_{t}(x,t) &  =   & y_{xx}(x,t), & {\rm in} \;
(0,1)\times(T_M,T),
\\y(0,t) & = & y(1,t) =  0,   &  {\rm in } \;\; (T_M,T),
\\ y(x,T_M) &  = & y_M(x,T_M),  & {\rm in } \;\; (0,1),
\end{array}
\right. \label{5.7}
\end{equation}
it is well known that
\begin{equation} \|\hat{y}_1(\cdot,T)\|_{L^2(0,1)} \leq \|y_M(\cdot,T_M)\|_{L^2(0,1)} \leq
\frac{\epsilon}{2}, \label{5.8}
\end{equation}
where the last inequality comes from Theorem \ref{T6}.

Thanks to the Theorem \ref{MP}-i), it holds that $y_j(x,t) \geq 0$ a.e. $(x,t) \in Q_T$ for $j = 1,\ldots,M$. In particular,
$y_M(x,T_M) \geq 0$ a.e. $x \in (0,1)$.

Assume that $y_d \in L^2_+ (0,1)$ and $\epsilon > 0$ are fixed. Applying now Theorem \ref{T5} for the corresponding system of
type $(S2)$ in the domain $(0,1)\times(T_M,T)$ with initial
datum $y_M(\cdot,T_M)$, we deduce the existence of $v_2 \in
L^\infty_+((0,1)\times (T_M,T))$ and $r_2 \in
PC[T_M,T]$ such that the solution $\hat{y}_2(x,t)$ of the
system
\begin{equation} \left\{
\begin{array}{lcll}
y_{t}(x,t) &  =   & y_{xx}(x,t) + v_2(x,t) y(x,t)
\chi_{\omega(t)}(x), & {\rm in} \; (0,1)\times(T_M,T),
\\y(0,t) & = & y(1,t) =  0,   &  {\rm in } \;\; (T_M,T),
\\ y(x,T_M) &  = & y_M(x,T_M),  & {\rm in } \;\; (0,1),
\end{array}
\right. \label{5.9}
\end{equation}
satisfies
\begin{equation}
\parallel \hat{y}_2 (\cdot, T) - \left(y_d +
\hat{y}_1(\cdot,T)\right)
\parallel_{L^2(0,1)} < \frac{\epsilon}{2}. \label{5.10}
\end{equation}

We finish the proof, by noticing that

\begin{equation} y_{y_0,v,r}(x,t) = \left\{
\begin{array}{ll}
\hat{y}_1(x,t), & \ (x,t) \in (0,1) \times (0,T_M],
\\ \hat{y}_2(x,t),    & \ (x,t) \in (0,1) \times [T_M,T),
\end{array}
\right. \label{5.11}
\end{equation}
with

\begin{equation} v(x,t) = \left\{
\begin{array}{ll}
-m_1, & \ (x,t) \in (0,1) \times (0,T_1),
\\ -m_2,    & \ (x,t) \in (0,1) \times (T_1,T_2),
\\ & \vdots
\\ -m_{j},    & \ (x,t) \in (0,1) \times (T_{j-1},T_j),
\\ & \vdots
\\ -m_M,    & \ (x,t) \in (0,1) \times (T_{M-1},T_M),
\\ v_2(x,t),    & \ (x,t) \in (0,1) \times (T_M,T),
\end{array}
\right. \label{5.12}
\end{equation}
and

\begin{equation} r(t) = \left\{
\begin{array}{ll}
0, & \ t \in (0,T_1),
\\ l,    & \ t \in (T_1, T_{2}),
\\ & \vdots
\\ jl,    & \ t \in (T_{j}, T_{j+1}),
\\ & \vdots
\\ (M-2)l,    & \ t \in (T_{M-2}, T_{M-1}),
\\ 1-l,    & \ t \in (T_{M-1},T_M),
\\ r_2(t), & \ t \in (T_M,T).
\end{array}
\right. \label{5.13}
\end{equation}

Clearly, $v \in L^\infty(Q_T)$ and $r \in PC[0,T]$. Moreover,
thanks to (\ref{5.8}) and (\ref{5.10}), it is satisfied
\[ \parallel y_{y_0,v,r}(\cdot, T) - y_d
\parallel_{L^2(0,1)} = \parallel y_{y_0,v,r}(\cdot, T) - (y_d + \hat{y}_1(\cdot,T)) + \hat{y}_1(\cdot,T)
\parallel_{L^2(0,1)} \leq \]
\[ \leq \parallel y_{y_0,v,r}(\cdot, T) - \left(y_d
+\hat{y}_1(\cdot,T)\right) \parallel_{L^2(0,1)}+ \|\hat{y}_1(\cdot,T)\|_{L^2(0,1)}
=  \]
\begin{equation} = \parallel \hat{y}_2(\cdot,T) - \left(y_d
+\hat{y}_1(\cdot,T)\right) \parallel_{L^2(0,1)}+ \|\hat{y}_1(\cdot,T)\|_{L^2(0,1)}
< \epsilon, \label{5.14}
\end{equation}
as we were looking for.\rule{2mm}{2mm}

\subsubsection{Auxiliary results.}

Our proof for Theorem \ref{T6} is quite technical and relies on the following results:

\begin{proposition}
Assume that $v \in C[0,1]$ with $v(x) \leq 0$ for all $x \in[0,1]$ and $y_0 \in C^{2+\sigma}[0,1]$, with $\sigma \in (0,1)$ verifies
$y_0(0)= y_0(1) = y_0''(0) = y_0''(1) = 0$ and $y_0(x) \geq 0$ for all $x \in[0,1]$.
Then, the classical solution of the problem
\begin{equation} \left\{
\begin{array}{lcll}
y_{t}(x,t) &  =   & y_{xx}(x,t) + v(x) y(x,t), & {\rm
in} \; Q_T,
\\y(0,t) & = & y(1,t) =  0,   &  {\rm in } \;\; (0,T),
\\ y(x,0) &  = & y_0 (x),  & {\rm in } \;\; (0,1),
\end{array}
\right. \label{5.20}
\end{equation}
satisfies the following properties:
\begin{itemize}
\item[a)] $y_t(x,t) \leq \max_{x \in[0,1]}(y_0''(x))^+, \ \ \ \forall \ (x,t) \in \overline{Q_T} $,
\item[b)] $0 \leq  y_x(0,t) \leq e \cdot \max_{x \in[0,1]} (y_0'(x)e^{y_0(x)}), \ \ \ \forall \ t \in [0,T]$ and
\item[c)] $e \cdot \min_{x \in[0,1]} (y_0'(x)e^{y_0(x)}) \leq  y_x(1,t) \leq 0, \ \ \ \forall \ t \in [0,T]$.
\end{itemize}
\label{P1}
\end{proposition}

{\bf Proof of Proposition \ref{P1}.} In the case $y_0 \equiv 0$ the result is clearly valid, because $y \equiv 0$. So, let us suppose $y_0 \not\equiv 0$.

The existence and uniqueness of a classical solution $y \in C^{2+\sigma,1+\sigma/2}(\overline{Q_T})$ for problem (\ref{5.20}) is a consequence of \cite[Theorem 5.2, p. 320]{14}. Moreover, by the Theorem \ref{MP}-i) and ii), it is also known that
\[0 \leq y(x,t) \leq \|y_0\|_{L^\infty(0,1)} \ \ \forall \ (x,t) \in \overline{Q_T}.\]

To prove $a)$, we differentiate (\ref{5.20}) with respect to the time variable
to deduce that $y_t$ is the solution of the problem
\begin{equation} \left\{\begin{array}{lcll}
z_{t}(x,t) &  =   & z_{xx}(x,t) + v(x) z(x,t), & {\rm
in} \; Q_T,
\\z(0,t) & = & z(1,t) =  0,   &  {\rm in } \;\; (0,T),
\\ z(x,0) &  = & y''_0 (x)+v(x)y_0(x),  & {\rm in } \;\; (0,1).
\end{array}
\right. \label{5.25}
\end{equation}
Using once more the Theorem \ref{MP}-ii) and taking into account that $v$ is nonpositive and $y_0$ is nonnegative in $[0,1]$,
we get
\[ y_t(x,t) \leq  \max_{x \in[0,1]}(y_0''(x)+v(x)y_0(x))^+ \leq \max_{x \in[0,1]}(y_0''(x))^+ \ \ \ \forall \ (x,t) \in \overline{Q_T}.\]

The lower bound of $b)$ is a direct consequence of the definition, the homogeneous boundary condition and that $y$ is nonnegative.
To derive the upper bound, we will use Bernstein's method (see \cite[p. 414 and 537]{14}), by introducing the auxiliary function
$w(x,t)= e^{y(x,t)} + \rho e^{1-x} -1$ with $\rho = \max_{x \in[0,1]}(y_0'(x)e^{y_0(x)})$. Let us remark that $\rho > 0$,
thanks to the hypotheses on $y_0$.

By direct calculation it can checked that $w \in C^{2,1}(\overline{Q_T})$ and satisfies the following PDE:
\[ w_t(x,t)-w_{xx}(x,t) = \left(v(x)y(x,t)-(y_x(x,t))^2\right)e^{y(x,t)} - \rho e^{1-x}.\]
Let us denote by $(x_0,t_0)$ a point in $\overline{Q_T}$ where $w$ takes its maximum value. We will see that $x_0 = 0$.
If $(x_0,t_0) \in Q_T$, it is well known that $w_t(x_0,t_0) = 0$ and $w_{xx}(x_0,t_0) \leq 0$, and using the previous PDE we arrive to the contradiction
\[ 0 \leq w_t(x_0,t_0)-w_{xx}(x_0,t_0) \leq - \rho e^{1-x_0} < 0.\]
The same argumentation can be used to exclude the case $t_0 = T$, because the unique difference is that $w_t(x_0,t_0) \geq 0$. Consequently, $(x_0,t_0)$ must be a point of the ``parabolic" boundary of $Q_T$, that is, $x_0 = 0$, $x_0 = 1$ or $t_0 = 0$. The last two possibilities can be excluded by noticing that
$w(0,t) = \rho e \geq \rho = w(1,t),$ $w(x,0) = e^{y_0(x)} + \rho e^{1-x} -1,$ and taking into account the choice for $\rho$:
\[ w_x(x,0) = y'_0(x)e^{y_0(x)} - \rho e^{1-x} \leq \rho - \rho e^{1-x} \leq 0, \ \ \forall x \in [0,1]. \]
Therefore, $w(x,0)$ is non-increasing in $[0,1]$ and $w(x,0) \leq w(0,0) = \rho e,$ for all $x \in [0,1]$.
Consequently, $w(x,t) \leq \rho e  = w(0,t),$ for all $(x,t) \in \overline{Q_T}$ and hence, $w_x(0,t) \leq 0$ for all $t \in [0,1]$.
This is equivalent to say that $y_x(0,t) \leq \rho e$ for all $t \in [0,1]$.

The proof for $c)$ is similar, using now the auxiliary function $\tilde{w}(x,t)= e^{y(x,t)} - \tilde{\rho} e^{x} -1$
with \[\tilde{\rho} = \min_{x \in[0,1]}(y_0'(x)e^{y_0(x)}) < 0.\rule{2mm}{2mm} \]

\begin{corollary}
Under the same hypotheses of Proposition \ref{P1}, except that $v \in L^\infty(0,1)$, the conclusions $a)-c)$ remain valid, for almost
every point.
\label{C1}
\end{corollary}

{\bf Proof of Corollary \ref{C1}.} By the classical convolution technique, it can be deduced the existence of a sequence $\{v_k\}_{k=1}^\infty \subset C[0,1]$ such that $-\|v\|_{L^\infty(0,1)} \leq v_k(x) \leq 0$ for all $k=1,2,\ldots,$ in $[0,1]$ and $v_k \rightarrow v$ in $L^q(0,1)$, as $k \rightarrow + \infty$, for some $q > 2$. If we denote by $y_{v_k} \in C^{2,1}(\overline{Q_T})$ the classical solution of problem (\ref{5.20}) with $v = v_k$, applying \cite[Theorem 9.1, p. 341]{14} to the difference $y_{v_k}-y$ and using that $0 \leq y_{v_k}(x,t) \leq \|y_0\|_{L^\infty(0,1)}$ in $Q_T$, by Theorem \ref{MP}, we get
\[ \|y_{v_k} -y\|_{W^{2,1}_q(Q_T)} \leq C \|v_k-v\|_{L^q(Q_T)} \rightarrow 0 \ \ \mbox{ as } \ k \rightarrow + \infty.\]
In particular,
\[ (y_{v_k})_t \longrightarrow  y_t \ \ \mbox{ in } L^q(Q_T), \ \mbox{ as } \ k \rightarrow + \infty.\]
Combining the estimate above with \cite[Lemma 3.4, p. 82]{14} for the values $m = s = 1$, $r = 0$ and $q > 2$, we deduce
\[ (y_{v_k})_x(0,t) \longrightarrow  y_x(0,t) \ \ \mbox{ and} \ \ (y_{v_k})_x(1,t) \longrightarrow  y_x(1,t)  \mbox{ in } L^q(0,T), \ \mbox{ as } \ k \rightarrow + \infty.\]
Using now the corresponding expressions $a)-c)$ for each $y_{v_k}$ and the fact that the lower and upper bounds do not depend on $v_k$, we can pass to the limit in them and arrive to the same ones for $y$ that will be satisfied for almost
every point. \rule{2mm}{2mm}

\subsubsection{Proof of Theorem \ref{T6}.}

The conclusion of the Theorem \ref{T6} is clearly valid if the $L^2-$norm of the initial datum is small enough compared with $\epsilon$. More precisely, if $\|y_0\|_{L^2(0,1)} \leq \frac{\epsilon}{2} e^{T\pi^2}$, it is sufficient to take $T_1 = T_2 = \ldots = T_M = 0$ and $m_1 = m_2 = \ldots = m_M = 0$. Of course, we want to deal with the general case.

Let us begin with a useful observation. As we have mentioned before, it is well known that for each $y_0 \in L^2_+(0,1)$, we can construct an approximating sequence $\{y_{0k}\}_{k=1}^\infty \subset C^{2+\sigma}[0,1]$ with compact support in $(0,1)$ such that $y_{0k}(x) \geq 0$ for all $x \in [0,1]$ and $k,$ with
$y_{0k} \rightarrow y_0$ in $L^2(0,1)$, as $k \rightarrow + \infty$. If we denote by $y_{y_0}$ the unique solution of (\ref{5.20})
and by $y_{y_{0k}}$ the unique solution of the same problem with initial datum $y_{0k}$ instead of $y_0$, subtracting both PDE, multiplying it by
the difference $y_{y_0}-y_{y_{0k}}$ and integrating by parts, we arrive to the expression
\[\frac{1}{2}\int_0^1 (y_{y_0}(x,T)-y_{y_{0k}}(x,T))^2 dx -\frac{1}{2}\int_0^1 (y_0(x)-y_{0k}(x))^2 dx = \]
\[ = - \int_{Q_T} ((y_{y_0}(x,t)-y_{y_{0k}}(x,t))_x)^2 dx dt + \int_{Q_T} v(x)(y_{y_0}(x,t)-y_{y_{0k}}(x,t))^2 dx dt \leq 0, \]
thanks to the hypothesis $v(x) \leq 0, \ $ a.e. $x \in (0,1)$. Therefore,
\[ \int_0^1 (y_{y_0}(x,T)-y_{y_{0k}}(x,T))^2 dx \leq \int_0^1 (y_0(x)-y_{0k}(x))^2 dx, \]
and, consequently, it is possible to make the $L^2-$norm of the difference of the solutions at time $T$ as small as needed, by taking
the $L^2-$norm of the initial conditions small enough. As a consequence of this observation and taking into account the conclusion of the Theorem \ref{T6} that we want to derive, we will assume in the sequel that the initial data of the problems appearing along the proof are nonnegative regular functions with compact support in $(0,1)$.

Now, we are in conditions to prove the existence of the nonnegative constants $m_1,\ldots, m_M$, the intermediate times $0 \leq T_1 \leq T_2 \leq \ldots \leq T_M < T$ and the corresponding solutions $y_j$ for $j = 1,\ldots,M$ (see (\ref{5.5})-(\ref{5.6-1})), such that the conclusion of the Theorem \ref{T6} holds.

For that purpose, we consider the family of problems
\begin{equation} \left\{
\begin{array}{lcll}
y_{t}(x,t) &  =   & y_{xx}(x,t) - m y(x,t) \chi_{(0,l)}(x), & {\rm
in} \; Q_T,
\\y(0,t) & = & y(1,t) =  0,   &  {\rm in } \;\; (0,T),
\\ y(x,0) &  = & y_0 (x),  & {\rm in } \;\; (0,1),
\end{array}
\right. \label{5.30}
\end{equation}
and denote its solution by $y_{m,1}$ for each $m \geq 0$. Multiplying the PDE by $y_{m,1}$ and integrating by parts in $Q_T$, we get as above
\[\frac{1}{2}\int_0^1 (y_{m,1}(x,T))^2 dx -\frac{1}{2}\int_0^1 (y_0(x))^2 dx = \]
\[ = - \int_{Q_T} ((y_{m,1})_x(x,t))^2 dx dt - m \int_0^T \int_0^l (y_{m,1}(x,t))^2 dx dt. \]
Consequently,
\begin{equation}
\int_0^T \int_0^l (y_{m,1}(x,t))^2 dx dt \leq \frac{1}{2m}\int_0^1 (y_0(x))^2 dx.
\label{5.40}
\end{equation}
From this inequality we deduce that the sequence $g_m(t) = \sqrt{\int_0^l (y_{m,1}(x,t))^2 dx}$
converges towards $0$ as $m \rightarrow + \infty$ in $L^2(0,T)$ and (taking a subsequence, if necessary) we
have that $g_m(t) \rightarrow 0$ as $m \rightarrow + \infty$ for a.e. $t \in (0,T)$. So, given $\epsilon > 0$, there exists
$T_1 \in (0,T)$ and $m_1 \gg 0$ such that $g_{m_1}(T_1) < \frac{\epsilon}{2\sqrt{2M-1}}$. Denoting $y_1 = y_{m_1,1}$ this is equivalent to say

\begin{equation}
\int_0^l (y_1(x,T_1))^2 dx \leq \frac{\epsilon^2}{4(2M-1)}.
\label{5.50}
\end{equation}

Following the same scheme, we consider now the family of problems
\begin{equation} \left\{
\begin{array}{lcll}
y_{t}(x,t) &  =   & y_{xx}(x,t) - m y(x,t) \chi_{(l,2l)}(x), & {\rm
in} \; (0,1)\times(T_1,T),
\\y(0,t) & = & y(1,t) =  0,   &  {\rm in } \;\; (T_1,T),
\\ y(x,T_1) &  = & y_{01}(x),  & {\rm in } \;\; (0,1),
\end{array}
\right. \label{5.60}
\end{equation}
with $y_{01}(x) = y_1(x,T_1)$, $y_1$ being the function obtained previously. If we denote its solution as $y_{m,2}$ for each $m \geq 0$,
repeating exactly the same argumentation, it is possible to derive the
existence of $T_2 \in (T_1,T)$ and $m_2 \gg 0$ such that $y_2 = y_{m_2,2}$ satisfies
\begin{equation}
\int_l^{2l} (y_2(x,T_2))^2 dx \leq \frac{\epsilon^2}{4(2M-1)}.
\label{5.70}
\end{equation}
In fact, we can select $T_2$ as close to $T_1$ as needed: more precisely, satisfying
\begin{equation} 0 < T_2 - T_1 \leq \min\{\frac{\epsilon^2}{8(2M-1)C_1}, T-T_1\}
\label{5.80}
\end{equation}
with
\begin{equation}C_1 = \|y_{01}\|_{L^\infty(0,1)}
\left(l \cdot \max_{x \in[0,1]}(y''_{01}(x))^+  + e \cdot \max_{x \in[0,1]} (y'_{01}(x)e^{y_{01}(x)})
\right).\label{5.90}
\end{equation}
At this point it is evident that we are selecting the initial datum $y_{01}$ to be a $C^{2+\sigma}$ function, with $\sigma \in (0,1)$,
with compact support in $(0,1)$ (see the beginning of this section). This is needed for the right-hand term of (\ref{5.90})
to be well defined. Let us remark that the terms appearing in (\ref{5.90}) are clearly related
with the conclusions of Proposition \ref{P1} (items $a)-b)$).

Under these conditions, we want to prove that

\begin{equation}
\int_0^{2l} (y_2(x,T_2))^2 dx \leq \frac{3\epsilon^2}{4(2M-1)}.
\label{5.100}
\end{equation}

Taking into account (\ref{5.70}), this will be true if we derive
\begin{equation}
\int_0^{l} (y_2(x,T_2))^2 dx \leq \frac{\epsilon^2}{2(2M-1)}.
\label{5.110}
\end{equation}
This is the most delicate point of the proof, where Corollary \ref{C1} is essential.
Again, we multiply the PDE of problem (\ref{5.60}) with $m = m_2$ by $y_2$ and integrate by parts in the domain $(x,t) \in (0,l) \times (T_1,T_2)$,
to obtain
\[\frac{1}{2}\int_0^l (y_2(x,T_2))^2 dx -\frac{1}{2}\int_0^l (y_{01}(x))^2 dx = \]
\begin{equation}= - \int_{T_1}^{T_2} \int_0^l((y_2)_x(x,t))^2 dx dt + \int_{T_1}^{T_2} (y_2)_x(l,t)y_2(l,t)dt, \label{5.120}
\end{equation}
where the boundary condition at $x=0$ has been utilized together with the fact that $(0,l) \cap (l,2l) = \emptyset$.
In particular, combining this relation with (\ref{5.50}), we have
\[ \int_0^l (y_2(x,T_2))^2 dx \leq \int_0^l (y_{01}(x))^2 dx + 2\int_{T_1}^{T_2} (y_2)_x(l,t)y_2(l,t)dt \leq \]
\begin{equation}
\leq \frac{\epsilon^2}{4(2M-1)}+ 2\int_{T_1}^{T_2} (y_2)_x(l,t)y_2(l,t)dt, \label{5.130}
\end{equation}
and the estimate (\ref{5.110}) will hold, if we are able to prove that
\begin{equation}
2\int_{T_1}^{T_2} (y_2)_x(l,t)y_2(l,t)dt \leq \frac{\epsilon^2}{4(2M-1)}. \label{5.140}
\end{equation}
Let us show that this holds thanks to the choice of $T_2$ (see (\ref{5.80})-(\ref{5.90})).
Using Corollary \ref{C1} (with $v(x) = - m_2 \cdot \chi_{(l,2l)}(x)$ that is only a bounded function) for problem (\ref{5.60}), it follows that
\[ (y_2)_x(l,t) = (y_2)_x(0,t) + \int_0^l (y_2)_{xx}(x,t) dx = (y_2)_x(0,t) + \int_0^l (y_2)_t(x,t) dx \leq \]
\[ \leq e \cdot \max_{x \in[0,1]} (y'_{01}(x)e^{y_{01}(x)}) + l \cdot \max_{x \in[0,1]} (y''_{01}(x))^+  \stackrel{def}{=} C_2 .\]
We can now combine this estimate with the fact that $0 \leq y_2 \leq \|y_{01}\|_{L^\infty(0,1)}$ thanks to the Theorem \ref{MP}-i) and ii), to get
\[ 2\int_{T_1}^{T_2} (y_2)_x(l,t)y_2(l,t)dt \leq 2\int_{T_1}^{T_2} C_2 \cdot \|y_{01}\|_{L^\infty(0,1)} dt = 2(T_2-T_1)C_1 \leq   \frac{\epsilon^2}{4(2M-1)},\]
that is exactly (\ref{5.140}).

Finally, we can repeat the argumentation to select the constant $m_j$, the time $T_j \in [T_{j-1},T)$ and the  solutions $y_j$, $j = 1,2,\ldots,M-1$,  verifying

\begin{equation}
\int_0^{jl} (y_j(x,T_j))^2 dx \leq \frac{(2j-1)\epsilon^2}{4(2M-1)},
\label{5.150}
\end{equation}
and arrive to
\begin{equation}
\int_0^1 (y_M(x,T_M))^2 dx \leq \frac{\epsilon^2}{4},
\label{5.160}
\end{equation}
as desired. \rule{2mm}{2mm}

\bibliographystyle{plain}

\end{document}